# Spill-Free Transfer and Stabilization of Viscous Liquid

**Iasson Karafyllis[*] and Miroslav Krstic[**]**

[*]Dept. of Mathematics, National Technical University of Athens, Zografou Campus, 15780, Athens, Greece, email: iasonkar@central.ntua.gr

[**]Dept. of Mechanical and Aerospace Eng., University of California, San Diego, La Jolla, CA 92093-0411, U.S.A., email: krstic@ucsd.edu

## Abstract

This paper studies the feedback stabilization problem of the motion of a tank that contains an incompressible, Newtonian, viscous liquid. The control input is the force applied on the tank and the overall system consists of two nonlinear Partial Differential Equations and two Ordinary Differential Equations. Moreover, a spill-free condition is required to hold. By applying the Control Lyapunov Functional methodology, a set of initial conditions (state space) is determined for which spill-free motion of the liquid is possible by applying an appropriate control input. Semi-global stabilization of the liquid and the tank by means of a simple feedback law is achieved, in the sense that for every closed subset of the state space, it is possible to find appropriate controller gains, so that every solution of the closed-loop system initiated from the given closed subset satisfies specific stability estimates. The closed-loop system exhibits an exponential convergence rate to the desired equilibrium point. The proposed stabilizing feedback law does not require measurement of the liquid level and velocity profiles inside the tank and simply requires measurements of: (i) the tank position error and tank velocity, (ii) the total momentum of the liquid, and (iii) the liquid levels at the tank walls. The obtained results allow an algorithmic solution of the problem of the spill-free movement and slosh-free settlement of a liquid in a vessel of limited height (such as water in a glass) by a robot to a pre-specified position, no matter how full the vessel is.

## 1. Introduction

One of the most important mathematical models that describes the free-surface flow of a thin layer of an incompressible, Newtonian liquid is the Saint-Venant system or shallow water model. The model has been used in many applications (river flows, tidal waves, oceans) and since its first derivation from first principles by Adhémar Jean Claude Barré de Saint-Venant in 1871 (see [1]), the model has been extended to take into account many types of forces acting on the body of the liquid other than gravity (e.g., viscous stresses, surface tension, friction forces; see [5,6,15,20,23,30]).

Feedback stabilization problems involving various variations of the Saint-Venant model have attracted the attention of many researchers during the last decades (see [2,3,4,7,8,9,10,11,12,21,24,25]). All works study the inviscid Saint-Venant model (i.e., the model that ignores viscous stresses and surface tension and takes into account gravity and friction forces), which provides a system of first-order hyperbolic Partial Differential Equations (PDEs). Many papers study stabilization problems for the linearization of this system around an equilibrium point, employing either the backstepping methodology (see [11,12]) or the Control Lyapunov Functional (CLF) methodology (providing local stabilization results for the original nonlinear system in many cases; see [2,3,4,8,9,10,21]).

The present work studies the problem of stabilization of the motion of a tank that contains an incompressible, Newtonian liquid. The problem is not new and has also been studied in [24] and we consider stabilization of both the state of the liquid (liquid level and velocity) and the state of the tank (tank position and velocity) by manipulating the acceleration of the tank. However, in contrast with [24], we do not consider the inviscid Saint-Venant model but the viscous Saint-Venant model (i.e., the model that ignores friction forces and surface tension and takes instead into account gravity forces and viscous stresses). From a mathematical point of view this makes a huge difference: the resulting nonlinear model consists of two Ordinary Differential Equations (ODEs) and two PDEs which are no longer first-order hyperbolic PDEs but a system of one hyperbolic first-order PDE and one parabolic PDE. From a physical point of view the replacement of the inviscid Saint-Venant model with friction by the viscous, frictionless Saint-Venant model is completely justified: viscous stresses are negligible compared to friction forces when the liquid velocity is relatively large (e.g. rivers) but in tanks where the liquid velocity is small, viscous stresses become dominant. The viscosity term is nonlinear and is similar to the viscosity term appearing in compressible fluid flow models (see [18,22,26,27]). In addition, we assume that the tank walls have a specific height and we require that the movement be spill-free and slosh-free (in the sense that the liquid is transferred without residual end point sloshing; see [13,32] for the description of the slosh-free motion problem). In other words, we require that the liquid level in the tank does not exceed the height of the tank walls. To our knowledge this is the first paper that studies the spill-free movement problem: all other works (see [24,25]) do not take into account the limited height of the tank walls.

We exploit the similarity with compressible fluid models in order to apply the CLF methodology for the viscous Saint-Venant system. The CLF methodology was first used for global stabilization of nonlinear parabolic PDEs in [19] and was subsequently studied in [16,17,18] (but see also [8] for the presentation of the CLF methodology in finite-dimensional and infinite-dimensional systems). Recently, it was used in [18] for the stabilization of a nonlinear Navier-Stokes model of compressible fluid. As in [18], the CLF plays also the role of a barrier function and guarantees a positive lower bound for the liquid level. However, here the CLF also plays the role of a barrier function for the additional requirement of spill-free movement by providing an upper bound for the liquid level. The CLF is constructed by combining the mechanical energy of the system and the use of a specific transformation that has been used in the literature of compressible fluids (see [18,22,26,27]). It is important to note here that the nonlinearity of the control problem hinders the application of standard feedback design methodologies like backstepping (see [28,29,31]).

The CLF allows the determination of a specific set $X$ of initial conditions for which spill-free movement of the tank is possible (by applying an appropriate input): this is the state space of our problem. We achieve semi-global stabilization of the liquid and the tank by means of a simple feedback law. The term "semi-global" here refers to the state space $X$ : for every closed subset of the state space $X$ we are in a position to find appropriate controller gains, so that every classical solution of the closed-loop system initiated from the given closed subset satisfies specific stability estimates. More specifically, an exponential convergence rate is achieved for the closed-loop system. The stabilizing feedback laws do not require measurement of the liquid level and velocity spatial profiles in the tank and simply require measurements of: (i) the tank position error and tank velocity, (ii) the total momentum of the liquid, and (iii) the liquid levels at the tank walls (boundary). To our knowledge, this is the first paper in the literature that achieves stabilization of the nonlinear viscous Saint-Venant system. In other words,

1) the nonlinear viscous Saint-Venant model is studied (instead of its linearization),
2) a simple feedback law (with minimum measurement requirements) is obtained that guarantees exponential convergence to the desired equilibrium point,

3) a spill-free movement of the liquid in the tank is guaranteed.

The structure of the paper is as follows. Section 2 is devoted to the presentation of the control problem. Section 3 contains the statements of the main result of the paper (Theorem 1). Moreover, in Section 3 we also provide insights for the proof of the main result: we explain the construction of the CLF, the selection of the state space and we provide some auxiliary results that are used in the proof. Section 4 provides an algorithmic solution of the problem of the spill-free movement and slosh-free settlement of water in a glass, by a robot, to a pre-specified glass position. We show that *no matter how full the glass is, the robot can transfer the glass without spilling out water and without residual end point sloshing.* The proofs of all results are provided in Section 5. Finally, Section 6 gives the concluding remarks of the present work.

**Notation.** Throughout this paper, we adopt the following notation.

* $\Re_+ = [0,+\infty)$ denotes the set of non-negative real numbers.

* Let $S \subseteq \Re^n$ be an open set and let $A \subseteq \Re^n$ be a set that satisfies $S \subseteq A \subseteq cl(S)$. By $C^0(A;\Omega)$, we denote the class of continuous functions on $A$, which take values in $\Omega \subseteq \Re^m$. By $C^k(A;\Omega)$, where $k \ge 1$ is an integer, we denote the class of functions on $A \subseteq \Re^n$, which takes values in $\Omega \subseteq \Re^m$ and has continuous derivatives of order $k$. In other words, the functions of class $C^k(A;\Omega)$ are the functions which have continuous derivatives of order $k$ in $S = \mathrm{int}(A)$ that can be continued continuously to all points in $\partial S \cap A$. When $\Omega = \Re$ then we write $C^0(A)$ or $C^k(A)$. When $I \subseteq \Re$ is an interval and $G \in C^1(I)$ is a function of a single variable, $G'(h)$ denotes the derivative with respect to $h \in I$.

* Let $I \subseteq \Re$ be an interval, let $a < b$ be given constants and let $u : I \times [a,b] \to \Re$ be a given function. We use the notation $u[t]$ to denote the profile at certain $t \in I$, i.e., $(u[t])(x) = u(t,x)$ for all $x \in [a,b]$. When $u(t,x)$ is (twice) differentiable with respect to $x \in [a,b]$, we use the notation $u_x(t,x)$ ($u_{xx}(t,x)$) for the (second) derivative of $u$ with respect to $x \in [a,b]$, i.e., $u_x(t,x) = \frac{\partial u}{\partial x}(t,x)$ ($u_{xx}(t,x) = \frac{\partial^2 u}{\partial x^2}(t,x)$). When $u(t,x)$ is differentiable with respect to $t$, we use the notation $u_t(t,x)$ for the derivative of $u$ with respect to $t$, i.e., $u_t(t,x) = \frac{\partial u}{\partial t}(t,x)$.

* Given a set $U \subseteq \Re^n$, $\chi_U$ denotes the characteristic function of $U$, i.e. the function defined by $\chi_U(x) := 1$ for all $x \in U$ and $\chi_U(x) := 0$ for all $x \notin U$. The sign function $\mathrm{sgn} : \Re \to \Re$ is the function defined by the relations $\mathrm{sgn}(x) = 1$ for $x > 0$, $\mathrm{sgn}(0) = 0$ and $\mathrm{sgn}(x) = -1$ for $x < 0$.

* Let $a < b$ be given constants. For $p \in [1,+\infty)$, $L^p(a,b)$ is the set of equivalence classes of Lebesgue measurable functions $u : (a,b) \to \Re$ with $\|u\|_p := \left( \int_a^b |u(x)|^p \, dx \right)^{1/p} < +\infty$. $L^\infty(a,b)$ is the set of equivalence classes of Lebesgue measurable functions $u : (a,b) \to \Re$ with $\|u\|_\infty := \underset{x \in (a,b)}{ess\sup}\left( |u(x)| \right) < +\infty$. For an integer $k \ge 1$, $H^k(a,b)$ denotes the Sobolev space of functions in $L^2(a,b)$ with all its weak derivatives up to order $k \ge 1$ in $L^2(a,b)$.

## 2. Description of the Problem

We consider a one-dimensional model for the motion of a tank. The tank contains a viscous, Newtonian, incompressible liquid. The tank is subject to a force that can be manipulated. We assume that the liquid pressure is hydrostatic and consequently, the liquid is modeled by the one-dimensional (1-D) viscous Saint-Venant equations, whereas the tank obeys Newton's second law.

The control objective is to drive asymptotically the tank to a specified position without liquid spilling out and having both the tank and the liquid within the tank at rest. Figure 1 shows a picture of the problem.

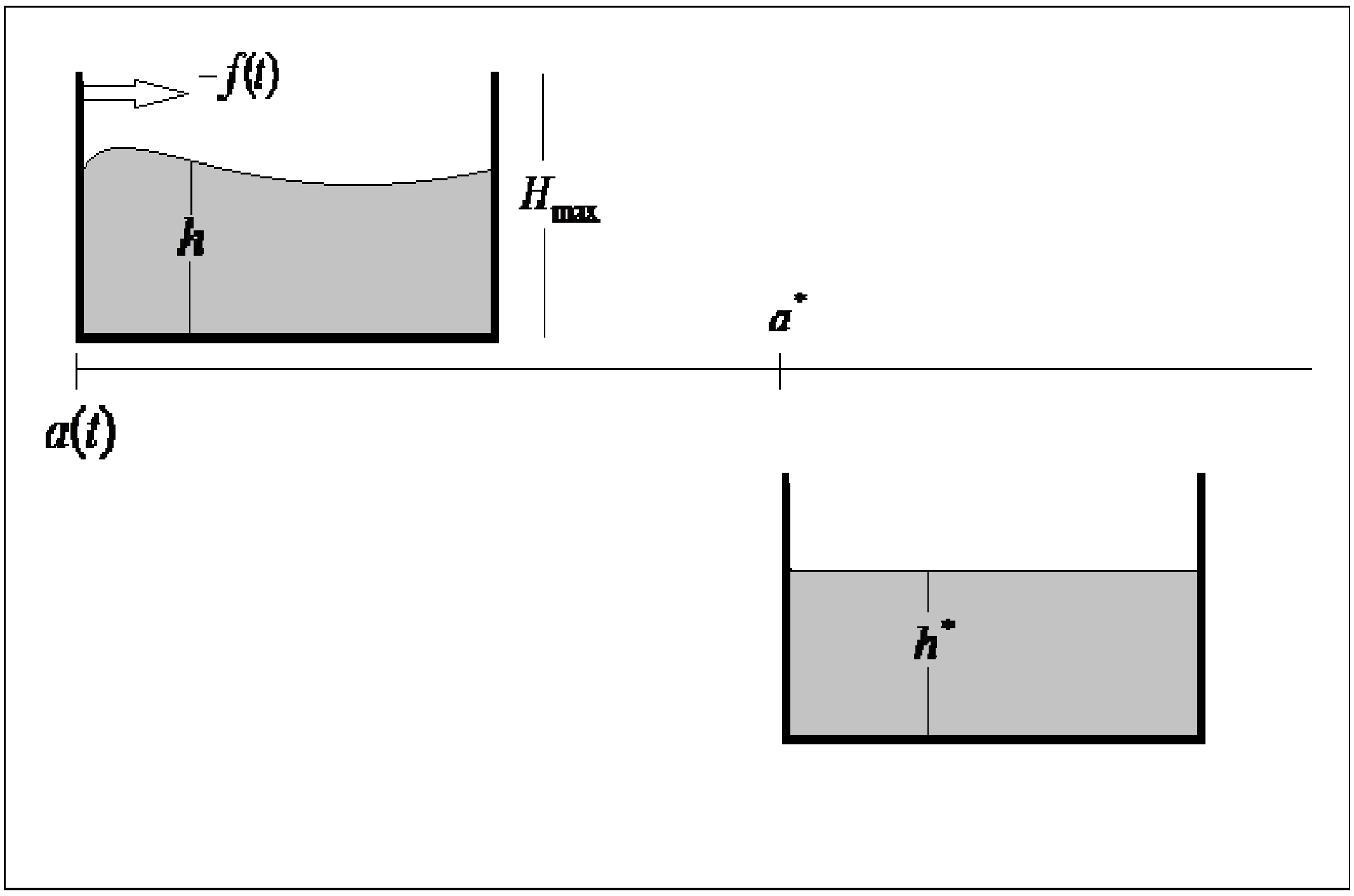

**Figure 1:** The control problem.

Let the position of the left side of the tank at time $t \geq 0$ be $a(t)$ and let the length of the tank be $L > 0$ (a constant). The equations describing the motion of the liquid within the tank are

$$H_t + (HV)_z = 0 \text{, for } t > 0,\ z \in \left[a(t), a(t)+L\right] \tag{1}$$

$$(HV)_t + \left(HV^2 + \frac{1}{2} g H^2\right)_z = \mu\left(HV_z\right)_z \text{, for } t > 0,\ z \in \left(a(t), a(t)+L\right) \tag{2}$$

where $H(t,z) > 0$, $V(t,z) \in \Re$ are the liquid level and the liquid velocity, respectively, at time $t \geq 0$ and position $z \in \left[a(t), a(t)+L\right]$, while $g, \mu > 0$ (constants) are the acceleration of gravity and the kinematic viscosity of the liquid, respectively.

The liquid velocities at the walls of the tank coincide with the tank velocity, i.e., we have:

$$V(t, a(t)) = V(t, a(t)+L) = w(t) \text{, for } t \geq 0 \tag{3}$$

where $w(t) = \dot{a}(t)$ is the velocity of the tank at time $t \geq 0$. Moreover, applying Newton's second law for the tank we get

$$\ddot{a}(t) = -f(t)\text{, for } t > 0 \tag{4}$$

where $-f(t)$, the control input to the problem, is the force exerted on the tank at time $t \geq 0$ divided by the total mass of the tank. The conditions for avoiding the liquid spilling out of the tank are:

$$\begin{aligned} &H(t, a(t)) < H_{\max} \\ &H(t, a(t) + L) < H_{\max} \end{aligned} \tag{5}$$

where $H_{\max} > 0$ is the height of the tank walls. The variable $f(t) \in \Re$ is the control input. Applying the transformation

$$\begin{aligned} &v(t, x) = V(t, a(t) + x) - w(t) \\ &h(t, x) = H(t, a(t) + x) \\ &\xi(t) = a(t) - a^* \end{aligned} \tag{6}$$

where $a^* \in \Re$ is the specified position (a constant) where we want to have the left side of the tank, we obtain the model

$$\dot{\xi} = w \quad , \quad \dot{w} = -f\text{, for } t \geq 0 \tag{7}$$

$$h_t + (hv)_x = 0\text{, for } t > 0,\ x \in [0, L] \tag{8}$$

$$(hv)_t + \left( hv^2 + \frac{1}{2} g h^2 \right)_x = \mu \left( h v_x \right)_x + hf\text{, for } t > 0,\ x \in (0, L) \tag{9}$$

$$v(t, 0) = v(t, L) = 0\text{, for } t \geq 0 \tag{10}$$

where the control input $f$ appears additively in the second equation of (7) and multiplicatively in (9). Moreover, the conditions (5) for avoiding the liquid spilling out of the tank become:

$$\max\left(h(t, 0), h(t, L)\right) < H_{\max}\text{, for } t \geq 0 \text{ --- Condition for no spilling out} \tag{11}$$

We consider classical solutions for the PDE-ODE system (7)-(10), i.e., we consider functions $\xi \in C^1\left(\Re_+\right) \cap C^2\left((0, +\infty)\right)$, $w \in C^0\left(\Re_+\right) \cap C^1\left((0, +\infty)\right)$, $v \in C^0\left([0, +\infty) \times [0, L]\right) \cap C^1\left((0, +\infty) \times [0, L]\right)$, $h \in C^1\left([0, +\infty) \times [0, L]; (0, +\infty)\right) \cap C^2\left((0, +\infty) \times (0, L)\right)$, with $v[t] \in C^2\left((0, L)\right)$ for each $t > 0$ that satisfy equations (7)-(10) for a given input $f \in C^0\left(\Re_+\right)$.

Using (8) and (10), we can prove that for every solution of (7)-(10) it holds that $\frac{d}{dt}\left( \int_0^L h(t, x)dx \right) = 0$ for all $t > 0$. Therefore, the total mass of the liquid $m > 0$ is constant. Therefore, without loss of generality, we assume that every solution of (7)-(10) satisfies the equation

$$\int_0^L h(t, x)dx \equiv m \tag{12}$$

The open-loop system (7)-(10), (12), i.e., system (7)-(10), (12) with $f(t) \equiv 0$, allows a continuum of equilibrium points, namely the points

$$h(x) \equiv h^*,\ v(x) \equiv 0, \text{ for } x \in [0, L] \tag{13}$$

$$\xi \in \Re,\ w = 0 \tag{14}$$

where $h^* = m / L$. We assume that the equilibrium points satisfy the condition for no spilling out (11), i.e., $h^* < H_{\max}$.

Our objective is to design a feedback law of the form

$$f(t) = F\left(h[t], v[t], \xi(t), w(t)\right), \text{ for } t > 0, \tag{15}$$

which achieves stabilization of the equilibrium point with $\xi = 0$. Moreover, we additionally require that the "spill-free condition" (11) holds for every $t \geq 0$.

The existence of a continuum of equilibrium points for the open-loop system given by (13), (14) means that the desired equilibrium point is not asymptotically stable for the open-loop system. Moreover, experience shows that there are smooth initial conditions $(\xi(0), w(0), h[0], v[0]) \in \Re \times \Re \times C^2\left([0, L]; (0, +\infty)\right) \times C^1\left([0, L]\right)$ with $(v[0])(0) = (v[0])(L) = 0$ for which it is not possible to avoid liquid spilling out of the tank-no matter what the applied input $f$ is. Thus, the described control problem is far from trivial.

## 3. Construction of the Feedback Law

### 3.1. The Control Lyapunov Functional (CLF)

Let $k, q > 0$ be position error and velocity gains (to be selected). Define the set $S \subset \Re^2 \times \left(C^0\left([0, L]\right)\right)^2$:

$$(\xi, w, h, v) \in S \Leftrightarrow \begin{cases} h \in C^0\left([0, L]; (0, +\infty)\right) \cap H^1(0, L) \\ v \in C^0\left([0, L]\right) \\ \int_0^L h(x)dx = m \\ (\xi, w) \in \Re^2, v(0) = v(L) = 0 \end{cases} \tag{16}$$

We define the following functionals for all $(\xi, w, h, v) \in S$:

$$V(\xi, w, h, v) := W(h, v) + E(h, v) + \frac{qk^2}{2}\xi^2 + \frac{q}{2}\left(w + k\xi\right)^2 \tag{17}$$

$$E(h, v) := \frac{1}{2}\int_0^L h(x)v^2(x)dx + \frac{1}{2}g\int_0^L \left(h(x) - h^*\right)^2 dx \tag{18}$$

$$W(h, v) := \frac{1}{2}\int_0^L h^{-1}(x)\left(h(x)v(x) + \mu h_x(x)\right)^2 dx + \frac{1}{2}g\int_0^L \left(h(x) - h^*\right)^2 dx \tag{19}$$

where $h^* = m / L$. We notice that:

- the functional $E$ is the mechanical energy of the liquid within the tank. Indeed, notice that $E$ is the sum of the potential energy ($\frac{1}{2}g\int_0^L \left(h(x)-h^*\right)^2 dx$) and the kinetic energy ($\frac{1}{2}\int_0^L h(x)v^2(x)dx$) of the liquid,

- the functional $W$ is a kind of mechanical energy of the liquid within the tank and has been used extensively in the literature of isentropic, compressible liquid flow (see [18,22,26,27]).

We intend to use the functional $V(\xi,w,h,v)$ defined by (17) as a CLF for the system. However, the functional $V(\xi,w,h,v)$ can also be used for the derivation of useful bounds for the function $h$. This is guaranteed by the following lemma.

**Lemma 1:** *Define the increasing* $C^1$ *function* $G:(0,+\infty)\to\left(-\frac{4}{3}h^*\sqrt{h^*},+\infty\right)$ *by means of the formula*

$$G(h):=\operatorname{sgn}\left(h-h^*\right)\left(\frac{2}{3}h\sqrt{h}-2h^*\sqrt{h}+\frac{4}{3}h^*\sqrt{h^*}\right) \tag{20}$$

*Let* $G^{-1}:\left(-\frac{4}{3}h^*\sqrt{h^*},+\infty\right)\to(0,+\infty)$ *be the inverse function of G and define*

$$c:=\frac{1}{\mu\sqrt{g}} \tag{21}$$

*Then for every* $(\xi,w,h,v)\in S$ *with* $V(\xi,w,h,v)<\frac{4}{3}\mu h^*\sqrt{gh^*}$*, the following inequality holds:*

$$G^{-1}\left(-cV(\xi,w,h,v)\right)\le h(x)\le G^{-1}\left(cV(\xi,w,h,v)\right),\text{ for all } x\in[0,L] \tag{22}$$

3.2. The state space

In order to be able to satisfy the condition for no spilling out (11) we need to restrict the state space. This becomes clear from the fact that the set $S$ contains states that violate the condition for no spilling out (11). Define

$$X:=\left\{(\xi,w,h,v)\in S : V(\xi,w,h,v)<R\right\} \tag{23}$$

where

$$R:=\frac{2\mu\sqrt{g}}{3}\left(2h^*\sqrt{h^*}+\sqrt{H_{\max}}\min\left(H_{\max}-3h^*,0\right)\right) \tag{24}$$

Notice that definitions (23), (24), the fact that $h^*<H_{\max}$ and Lemma 1 imply that for all $(\xi,w,h,v)\in X$ it holds that

$$0<G^{-1}\left(-cV(\xi,w,h,v)\right)\le h(x)\le G^{-1}\left(cV(\xi,w,h,v)\right)<H_{\max},\text{ for all } x\in[0,L] \tag{25}$$

Therefore, the condition for no spilling out (11) is automatically satisfied when $(\xi,w,h,v)\in X$. However, it should be noticed that requiring the state to be in $X$ in order to satisfy the condition for no spilling out (11) is a conservative approach: we actually require that the liquid level $h(x)$ is below the height of the walls $H_{\max}$ everywhere in the tank (i.e., for all $x \in [0, L]$; recall (25)) and not simply at the walls (i.e., at $x = 0$ and $x = 1$) as the condition for no spilling out (11) requires.

The state space of system (7)-(10), (12) is considered to be the set $X$. More specifically, we consider as state space the set $X \subset \Re^2 \times H^1(0,L) \times L^2(0,L)$ endowed with the norm of the underlying normed linear space $\Re^2 \times H^1(0,L) \times L^2(0,L)$, i.e.,

$$\|(\xi, w, h, v)\|_X = \left(\xi^2 + w^2 + \|h\|_2^2 + \|h'\|_2^2 + \|v\|_2^2\right)^{1/2} \tag{26}$$

## 3.3. Main result

Now we are in a position to state the first main result of the present work.

**Theorem 1:** *For every $r \in [0, R)$, where $R > 0$ is the constant defined by (24), and for every $\sigma, q, k > 0$ with*

$$k < q\theta \frac{G^{-1}(-cr)}{b + G^{-1}(-cr)} \tag{27}$$

*where $c > 0$ is defined by (21), $\theta := \dfrac{\sigma g}{g + \mu\sigma L}$, $b := \dfrac{4mL^2 H_{\max}}{\mu\pi^2}\theta$, there exist constants $M, \lambda > 0$ with the following property:*

**(P)** *Every classical solution of the PDE-ODE system (7)-(10), (12) and*

$$f(t) = -\sigma\left(2\int_0^L h(t,x)v(t,x)dx + \mu\left(h(t,L) - h(t,0)\right) - q\left(w(t) + k\xi(t)\right)\right), \text{ for } t > 0 \tag{28}$$

*with $V(\xi(0), w(0), h[0], v[0]) \le r$, satisfies $(\xi(t), w(t), h[t], v[t]) \in X$ and the following estimate for all $t \ge 0$:*

$$\left\|(\xi(t), w(t), h[t] - h^*\chi_{[0,L]}, v[t])\right\|_X \le M \exp(-\lambda t)\left\|(\xi(0), w(0), h[0] - h^*\chi_{[0,L]}, v[0])\right\|_X \tag{29}$$

Theorem 1 guarantees semi-global exponential stabilization of the state in the norm of $\Re^2 \times H^1(0,L) \times L^2(0,L)$ by means of the nonlinear feedback law (28). Here the term "semi-global" refers to the set $X$ for which spill-free transfer can be guaranteed: the constant $r < R$ can be arbitrarily close to $R$. It should be noticed that the feedback law (28) does not require the measurement of the whole liquid level and liquid velocity profile and requires the measurement of only four quantities:

- the tank position $\xi(t)$ and the tank velocity $w(t)$,
- the total liquid momentum $\int_0^L h(t,x)v(t,x)dx$, and
- the liquid level difference at the tank walls $h(t,L) - h(t,0)$.

The proof of Theorem 1 shows that as $r \to R$ we have $k \to 0$, $\lambda \to 0$ and $M \to +\infty$. This is expected because as $r \to R$ we may have initial conditions satisfying $V(\xi(0), w(0), h[0], v[0]) \le r$ for which the liquid in the container is agitated. Therefore, in such cases the movement of the container may cause additional agitation of the liquid which can make the liquid spill-out unavoidable.

### 3.4. Auxiliary results

For the proof of Theorem 1, we need some auxiliary lemmas. The first auxiliary lemma provides formulas for the time derivatives of the energy functionals defined by (18), (19).

**Lemma 2:** *For every classical solution of the PDE-ODE system (7)-(10), (12) the following equations hold for all* $t>0$*:*

$$\frac{d}{dt}E\left(h[t],v[t]\right)=-\mu\int_0^L h(t,x)v_x^2(t,x)dx+f(t)\int_0^L h(t,x)v(t,x)dx \tag{30}$$

$$\frac{d}{dt}W\left(h[t],v[t]\right)=-\mu g\int_0^L h_x^2(t,x)dx+f(t)\int_0^L\left(h(t,x)v(t,x)+\mu h_x(t,x)\right)dx \tag{31}$$

*where* $E,W$ *are the functionals defined by (18), (19), respectively.*

The two next auxiliary lemmas provide useful inequalities for the CLF $V$ defined by (17).

**Lemma 3:** *Let* $q,k>0$ *be given. Then there exists a non-decreasing function* $\Gamma:[0,R)\to(0,+\infty)$*, where* $R>0$ *is defined by (24), such that for every* $(\xi,w,h,v)\in X$ *with* $v\in H^1(0,L)$*, the following inequality holds:*

$$V(\xi,w,h,v)\le\Gamma(V(\xi,w,h,v))\left(\int_0^L h_x^2(x)dx+\int_0^L h(x)v_x^2(x)dx+\xi^2+\left(w+k\xi\right)^2\right) \tag{32}$$

**Lemma 4:** *Let* $q,k>0$ *be given. Then there exist non-decreasing functions* $G_i:[0,R)\to(0,+\infty)$*,* $i=1,2$*, where* $R>0$ *is defined by (24), such that for every* $(\xi,w,h,v)\in X$ *, the following inequality holds:*

$$\frac{V(\xi,w,h,v)}{G_2\left(V(\xi,w,h,v)\right)}\le\left\|(\xi,w,h-h^*\chi_{[0,L]},v)\right\|_X^2\le V(\xi,w,h,v)G_1\left(V(\xi,w,h,v)\right) \tag{33}$$

*where* $\|\cdot\|_X$ *is defined by (26).*

Inequality (32) provides an estimate of the dissipation rate of the Lyapunov functional for the closed-loop system (7)-(10), (12) with (28). On the other hand, inequalities (33) provide estimates of the Lyapunov functional in terms of the norm of the state space.

## 4. Can a Robot Move a Glass of Water Without Spilling Out Water?

We next consider the problem of the movement of a glass of water by means of a robotic arm, from near-rest to near-rest in finite (but not prespecified) time. The glass of water starting from an almost at rest state (both the glass and the water in the glass), is formulated as

$$\left\| (0, w(0), h[0] - h^* \chi_{[0,L]}, v[0]) \right\|_X \leq \varepsilon \tag{34}$$

where $\varepsilon > 0$ is a small number (tolerance) and $\| \cdot \|_X$ is the norm defined by (26). The robotic arm should move the glass to a pre-specified position without spilling out water of the glass and having the water in the glass almost still at the final time, i.e., at the final time $T > 0$ we must have

$$\left\| (\xi(T), w(T), h[T] - h^* \chi_{[0,L]}, v[T]) \right\|_X \leq \varepsilon \tag{35}$$

In other words, we require the spill-free and slosh-free motion of the glass. The initial condition is $\left( \xi(0), w(0), h[0], v[0] \right) \in S$, where $\xi(0) \neq 0$ (recall the definition of $S$ (16)). The problem can be solved by using Theorem 1 when the tolerance $\varepsilon > 0$ is sufficiently small. More specifically, we require that $\varepsilon > 0$ is small so that

$$\varepsilon^2 \max\left( \mu^2 \left( h^* - \varepsilon\sqrt{L} \right)^{-1}, g, \frac{3H_{\max}}{2} \right) < R \text{ and } \varepsilon < \frac{\min\left( h^*, H_{\max} - h^* \right)}{\sqrt{L}} \tag{36}$$

where $R := \frac{2\mu\sqrt{g}}{3}\left( 2h^*\sqrt{h^*} + \sqrt{H_{\max}} \min\left( H_{\max} - 3h^*, 0 \right) \right)$.

If inequalities (36) hold then we can follow the following algorithm:

<u>Step 1:</u> Pick numbers $r \in (0, R)$, $\sigma, q > 0$ with $q \leq \max\left( \mu^2 \left( h^* - \varepsilon\sqrt{L} \right)^{-1}, g, \frac{3H_{\max}}{2} \right)$ and $\varepsilon^2 \max\left( \mu^2 \left( h^* - \varepsilon\sqrt{L} \right)^{-1}, g, \frac{3H_{\max}}{2} \right) < r$ (this is possible due to (36)).

<u>Step 2:</u> Select $k > 0$ so that (27) holds and so that

$$k \leq \sqrt{\frac{2}{3q}} \frac{\sqrt{r - \varepsilon^2 \max\left( \mu^2 \left( h^* - \varepsilon\sqrt{L} \right)^{-1}, g, \frac{3H_{\max}}{2}, q \right)}}{|\xi(0)|} \tag{37}$$

Inequality (37) guarantees the inequality $V(\xi(0), w(0), h[0], v[0]) \leq r$. Indeed, this fact is a direct consequence of (34), (36), (37) and the following proposition.

**Proposition 1:** *Let $q, k > 0$ be given. Then for every $(\xi, w, h, v) \in S$ satisfying the inequality $\left\| (0, w, h - h^* \chi_{[0,L]}, v) \right\|_X \leq \varepsilon$ for some $\varepsilon > 0$ with $\varepsilon < \frac{\min\left( h^*, H_{\max} - h^* \right)}{\sqrt{L}}$, the following inequality holds:*

$$V(\xi,w,h,v)\le \max\left(\mu^2\left(h^*-\varepsilon\sqrt{L}\right)^{-1},g,\frac{3H_{\max}}{2},q\right)\left\|(0,w,h-h^*\chi_{[0,L]},v)\right\|_X^2+\frac{3qk^2}{2}\xi^2 \tag{38}$$

*where* $\|\cdot\|_X$ *is defined by (26).*

Let $M,\lambda>0$ be the constants involved in (29) and correspond to the selected parameters $r\in(0,R)$, $\sigma,q>0$ and $k>0$.

Step 3: Set $T=\frac{1}{\lambda}\ln\left(\frac{M\left|\xi(0)\right|+M\varepsilon}{\varepsilon}\right)$ and apply the feedback law (28) for $t\in(0,T]$. Inequality (29) implies that estimate (35) holds.

The application of the above algorithm can guarantee that the robotic arm will move the glass to the pre-specified position *without spilling out water of the glass and without residual end point sloshing no matter how small the difference* $H_{\max}-h^*>0$ *is*. Of course, the (minimization of the) final time $T>0$ depends on the (appropriate) selection of the parameters $r\in(0,R)$, $\sigma,q>0$ and $k>0$.

## 5. Proofs

We start by providing the proof of Lemma 1.

**Proof of Lemma 1:** The fact that the function $G:(0,+\infty)\to\left(-\frac{4}{3}h^*\sqrt{h^*},+\infty\right)$ defined by (20) is a $C^1$ increasing function follows from the fact that the equation $G(h)=\int_{h^*}^{h}\frac{\left|r-h^*\right|}{\sqrt{r}}dr$ holds for all $h>0$.

Let $(\xi,w,h,v)\in S$ with $V(\xi,w,h,v)<\frac{4}{3}\mu h^*\sqrt{gh^*}$ be given. Using the inequality $(hv+\mu h_x)^2\ge\frac{1}{2}\mu^2h_x^2-h^2v^2$ and definition (19), we obtain the following estimate:

$$W(h,v)\ge\frac{\mu^2}{4}\int_0^L h^{-1}(x)h_x^2(x)dx-\frac{1}{2}\int_0^L h(x)v^2(x)dx+\frac{1}{2}g\int_0^L\left(h(x)-h^*\right)^2dx \tag{39}$$

Using definitions (17), (18) and estimate (39), we obtain:

$$\frac{\mu^2}{4}\int_0^L h^{-1}(x)h_x^2(x)dx+g\int_0^L\left(h(x)-h^*\right)^2dx\le W(h,v)+E(h,v)\le V(\xi,w,h,v) \tag{40}$$

Let arbitrary $x,y\in[0,L]$ be given. Using the Cauchy-Schwarz inequality and the fact that $G'(h)=\frac{\left|h-h^*\right|}{\sqrt{h}}$ for all $h>0$ (a consequence of definition (20)), we get:

$$\left|G(h(x))-G(h(y))\right| \le \left|\int_y^x G'(h(s))h_s(s)ds\right| \le \int_{\min(x,y)}^{\max(x,y)} \left|G'(h(s))h_s(s)\right|ds$$
$$\le \left(\int_{\min(x,y)}^{\max(x,y)} h^{-1}(s)h_s^2(s)ds\right)^{1/2} \left(\int_{\min(x,y)}^{\max(x,y)} h(s)\left(G'(h(s))\right)^2 ds\right)^{1/2}$$
$$\le \left(\int_0^L h^{-1}(s)h_s^2(s)ds\right)^{1/2} \left(\int_0^L h(s)\left(G'(h(s))\right)^2 ds\right)^{1/2} \tag{41}$$
$$= \left(\int_0^L h^{-1}(s)h_s^2(s)ds\right)^{1/2} \left(\int_0^L \left|h(s)-h^*\right|^2 ds\right)^{1/2}$$

Since $\max\left\{ ab : \frac{\mu^2}{4}a^2 + gb^2 \le V \right\} = \frac{V}{\mu\sqrt{g}}$ for all $V \ge 0$, we obtain from (40), (41) and definition (21):

$$\left|G(h(x))-G(h(y))\right| \le cV(\xi,w,h,v), \text{ for all } x,y \in [0,L] \tag{42}$$

Since $\int_0^L h(x)dx = m$ (recall definition (16)) and since $h^* = m/L$, it follows from continuity of $h$ and the mean value theorem that there exists $x^* \in [0,L]$ such that $h(x^*) = h^*$. Moreover, since $G(h^*) = 0$ (a consequence of definition (20)), we get from (42) (with $y = x^*$):

$$-cV(\xi,w,h,v) \le G(h(x)) \le cV(\xi,w,h,v), \text{ for all } x \in [0,L] \tag{43}$$

Inequality (22) is a direct consequence of (43), the fact that the inverse function of $G$ is increasing and the fact that $cV\left(\xi,w,h,v\right) < \frac{4}{3}h^*\sqrt{h^*}$. The proof is complete. ⊲

We continue with the proof of Lemma 2.

**Proof of Lemma 2:** Using (18) we obtain for all $t > 0$:

$$\frac{d}{dt}E(h[t],v[t]) = \frac{1}{2}\int_0^L h_t(t,x)v^2(t,x)dx$$
$$+\int_0^L h(t,x)v(t,x)v_t(t,x)dx + g\int_0^L \left(h(t,x)-h^*\right)h_t(t,x)dx \tag{44}$$

Using (8), (9) we obtain:

$$v_t + vv_x + gh_x = \mu h^{-1}\left(hv_x\right)_x + f, \text{ for } t>0,\ x \in \left(0,L\right) \tag{45}$$

Combining (8), (44) and (45) we get for all $t > 0$:

$$\frac{d}{dt}E(h[t],v[t]) = -\frac{1}{2}\int_0^L \left(h(t,x)v(t,x)\right)_x v^2(t,x)dx + \mu\int_0^L v(t,x)\left(h(t,x)v_x(t,x)\right)_x dx$$
$$+f(t)\int_0^L h(t,x)v(t,x)dx - \int_0^L h(t,x)v^2(t,x)v_x(t,x)dx \tag{46}$$
$$-g\int_0^L h(t,x)v(t,x)h_x(t,x)dx - g\int_0^L \left(h(t,x)-h^*\right)\left(h(t,x)v(t,x)\right)_x dx$$

Integrating by parts and using (10), we get for all $t>0$:

$$\int_0^L v(t,x)\left(h(t,x)v_x(t,x)\right)_x dx = -\int_0^L h(t,x)v_x^2(t,x)dx$$
$$\int_0^L \left(h(t,x)-h^*\right)\left(h(t,x)v(t,x)\right)_x dx = -\int_0^L h(t,x)v(t,x)h_x(t,x)dx \tag{47}$$
$$\int_0^L \left(h(t,x)v(t,x)\right)_x v^2(t,x)dx = -2\int_0^L h(t,x)v^2(t,x)v_x(t,x)dx$$

Equation (30) is a direct consequence of (46) and (47).

Next define for all $t \geq 0$ and $x \in [0,L]$:

$$\varphi(t,x) = h(t,x)v(t,x) + \mu h_x(t,x) \tag{48}$$

Using (8), (9) and definition (48) we conclude that the following equation holds for all $t>0$ and $x \in (0,L)$:

$$\varphi_t(t,x) = h(t,x)f(t) - \left(\varphi(t,x)v(t,x) + \frac{1}{2}gh^2(t,x)\right)_x \tag{49}$$

Using (19) and definition (48), we obtain for all $t>0$:

$$\frac{d}{dt}W(h[t],v[t]) = -\frac{1}{2}\int_0^L h^{-2}(t,x)h_t(t,x)\varphi^2(t,x)dx$$
$$+\int_0^L h^{-1}(t,x)\varphi(t,x)\varphi_t(t,x)dx + g\int_0^L \left(h(x)-h^*\right)h_t(t,x)dx \tag{50}$$

Combining (B12) and (8), (49), we get for all $t>0$:

$$\frac{d}{dt}W(h[t],v[t]) = \frac{1}{2}\int_0^L h^{-2}(t,x)\varphi^2(t,x)\left(h(t,x)v(t,x)\right)_x dx$$
$$+f(t)\int_0^L \varphi(t,x)dx - \int_0^L h^{-1}(t,x)\varphi(t,x)\left(\varphi(t,x)v(t,x) + \frac{1}{2}gh^2(t,x)\right)_x dx \tag{51}$$
$$-g\int_0^L \left(h(x)-h^*\right)\left(h(t,x)v(t,x)\right)_x dx$$

Integrating by parts and using (10), we get for all $t>0$:

$$\int_0^L h^{-2}(t,x)\varphi^2(t,x)\left(h(t,x)v(t,x)\right)_x dx$$
$$= 2\int_0^L h^{-2}(t,x)h_x(t,x)\varphi^2(t,x)v(t,x)dx - 2\int_0^L h^{-1}(t,x)\varphi(t,x)\varphi_x(t,x)v(t,x)dx \qquad (52)$$

Combining (51), (52) and (47), we obtain for all $t>0$:

$$\begin{aligned}\frac{d}{dt}W(h[t],v[t]) &= \int_0^L h^{-2}(t,x)h_x(t,x)\varphi^2(t,x)v(t,x)dx \\ &-2\int_0^L h^{-1}(t,x)\varphi(t,x)\varphi_x(t,x)v(t,x)dx + f(t)\int_0^L \varphi(t,x)dx \\ &-\int_0^L h^{-1}(t,x)\varphi^2(t,x)v_x(t,x)dx \\ &+g\int_0^L \left(h(t,x)v(t,x)-\varphi(t,x)\right)h_x(t,x)dx\end{aligned} \qquad (53)$$

Equation (53) in conjunction with definition (48) gives for all $t>0$:

$$\begin{aligned}\frac{d}{dt}W(h[t],v[t]) &= -\int_0^L \left(\varphi^2(t,x)h^{-1}(t,x)v(t,x)\right)_x dx \\ &+ f(t)\int_0^L \left(h(t,x)v(t,x)+\mu h_x(t,x)\right)dx - g\mu\int_0^L h_x^2(t,x)dx\end{aligned} \qquad (54)$$

Equation (31) is a direct consequence of (54) and (10). The proof is complete. ⊲

Next, we provide the proofs of Lemma 3 and Lemma 4.

**Proof of Lemma 3:** Let arbitrary $(\xi,w,h,v)\in X$ with $v\in H^1(0,L)$ be given. Since $\int_0^L h(x)dx = m$ (recall definitions (16), (23)) and since $h^* = m/L$, it follows from continuity of $h$ and the mean value theorem that there exists $x^*\in[0,L]$ such that $h(x^*)=h^*$. Using the Cauchy-Schwarz inequality, we get for all $x\in[0,L]$:

$$h(x)-h^* = \int_{x^*}^{x} h_s(s)ds \Rightarrow \left|h(x)-h^*\right| \le \int_{\min(x,x^*)}^{\max(x,x^*)} \left|h_s(s)\right|ds \le \int_0^L \left|h_s(s)\right|ds \le \sqrt{L}\left(\int_0^L h_s^2(s)ds\right)^{1/2} \qquad (55)$$

Therefore, using (55), we obtain:

$$\int_0^L \left(h(x)-h^*\right)^2 dx \le L^2\int_0^L h_x^2(x)dx \qquad (56)$$

Using the inequality $(h(x)v(x)+\mu h_x(x))^2 \le 2h^2(x)v^2(x)+2\mu^2 h_x^2(x)$ and (25), we obtain:

$$\int_0^L h^{-1}(x)(h(x)v(x)+\mu h_x(x))^2 dx \le 2\int_0^L h(x)v^2(x)dx + 2\mu^2\int_0^L h^{-1}(x)h_x^2(x)dx$$
$$\le 2\int_0^L h(x)v^2(x)dx + \frac{2\mu^2}{G^{-1}\left(-cV(\xi,w,h,v)\right)}\int_0^L h_x^2(x)dx \tag{57}$$

where $G$ is defined by (20). Since $v(0)=v(L)=0$ (recall definitions (16), (23)), by virtue of Wirtinger's inequality we have:

$$\int_0^L v^2(x)dx \le \frac{L^2}{\pi^2}\int_0^L v_x^2(x)dx \tag{58}$$

Combining (58) and (25) we get:

$$\int_0^L h(x)v^2(x)dx \le G^{-1}\left(cV(\xi,w,h,v)\right)\int_0^L v^2(x)dx$$
$$\le \frac{L^2}{\pi^2}G^{-1}\left(cV(\xi,w,h,v)\right)\int_0^L v_x^2(x)dx \le \frac{L^2G^{-1}\left(cV(\xi,w,h,v)\right)}{\pi^2G^{-1}\left(-cV(\xi,w,h,v)\right)}\int_0^L h(x)v_x^2(x)dx \tag{59}$$

Using (17), (18), (19) and (56), (57), (59), we obtain:

$$V\left(\xi,w,h,v\right) = \frac{1}{2}\int_0^L h^{-1}(x)(h(x)v(x)+\mu h_x(x))^2 dx + g\int_0^L \left(h(x)-h^*\right)^2 dx$$
$$+\frac{1}{2}\int_0^L h(x)v^2(x)dx + \frac{qk^2}{2}\xi^2 + \frac{q}{2}\left(w+k\xi\right)^2$$
$$\le \frac{3}{2}\int_0^L h(x)v^2(x)dx + \left(\frac{\mu^2}{G^{-1}\left(-cV(\xi,w,h,v)\right)} + gL^2\right)\int_0^L h_x^2(x)dx \tag{60}$$
$$+\frac{qk^2}{2}\xi^2 + \frac{q}{2}\left(w+k\xi\right)^2$$
$$\le \Gamma\left(V(\xi,w,h,v)\right)\left(\int_0^L h(x)v_x^2(x)dx + \int_0^L h_x^2(x)dx + \xi^2 + \left(w+k\xi\right)^2\right)$$

where

$$\Gamma\left(s\right) := \max\left(\frac{3L^2G^{-1}\left(cs\right)}{2\pi^2G^{-1}\left(-cs\right)}, \frac{\mu^2}{G^{-1}\left(-cs\right)} + gL^2, \frac{qk^2}{2}, \frac{q}{2}\right), \text{ for } s\in[0,R) \tag{61}$$

The proof is complete. ◁

**Proof of Lemma 4:** Let arbitrary $(\xi,w,h,v)\in X$ be given. Using definitions (17), (18), (19), the inequalities $(h(x)v(x)+\mu h_x(x))^2 \le 2h^2(x)v^2(x)+2\mu^2h_x^2(x)$, $\left(w+k\xi\right)^2 \le 2w^2+2k^2\xi^2$ and (25), we obtain:

$$
\begin{aligned}
V\left(\xi,w,h,v\right) &= \frac{1}{2}\int_0^L h^{-1}(x)(h(x)v(x)+\mu h_x(x))^2 dx + g\int_0^L \left(h(x)-h^*\right)^2 dx \\
&+\frac{1}{2}\int_0^L h(x)v^2(x)dx + \frac{qk^2}{2}\xi^2 + \frac{q}{2}\left(w+k\xi\right)^2 \\
&\le \frac{3}{2}\int_0^L h(x)v^2(x)dx + \mu^2\int_0^L h^{-1}(x)h_x^2(x)dx + g\int_0^L \left(h(x)-h^*\right)^2 dx + \frac{3qk^2}{2}\xi^2 + qw^2 \qquad (62)\\
&\le \frac{3}{2}H_{\max}\int_0^L v^2(x)dx + \frac{\mu^2}{G^{-1}\left(-cV(\xi,w,h,v)\right)}\int_0^L h_x^2(x)dx + g\int_0^L \left(h(x)-h^*\right)^2 dx \\
&+\frac{3qk^2}{2}\xi^2 + qw^2
\end{aligned}
$$

Inequality (62) implies the left inequality (33) with $G_2\left(s\right) := \max\left(\frac{3}{2}H_{\max}, \frac{\mu^2}{G^{-1}\left(-cs\right)}, \frac{3qk^2}{2}, q, g\right)$ for $s \in [0,R)$.

Using definitions (17), (18), (19), the inequalities $v(x)h_x(x) \ge -\frac{3}{4\mu}h(x)v^2(x) - \frac{\mu}{3}h^{-1}(x)h_x^2(x)$, $w\xi \ge -\frac{3k}{4}\xi^2 - \frac{1}{3k}w^2$ and (25), we obtain:

$$
\begin{aligned}
V\left(\xi,w,h,v\right) &= \frac{1}{2}\int_0^L h^{-1}(x)(h(x)v(x)+\mu h_x(x))^2 dx + g\int_0^L \left(h(x)-h^*\right)^2 dx \\
&+\frac{1}{2}\int_0^L h(x)v^2(x)dx + \frac{qk^2}{2}\xi^2 + \frac{q}{2}\left(w+k\xi\right)^2 \\
&\ge \frac{1}{4}\int_0^L h(x)v^2(x)dx + \frac{\mu^2}{6}\int_0^L h^{-1}(x)h_x^2(x)dx + g\int_0^L \left(h(x)-h^*\right)^2 dx + \frac{qk^2}{4}\xi^2 + \frac{q}{6}w^2 \qquad (63)\\
&\ge \frac{1}{4}G^{-1}\left(-cV\left(\xi,w,h,v\right)\right)\int_0^L v^2(x)dx + \frac{\mu^2}{6H_{\max}}\int_0^L h_x^2(x)dx + g\int_0^L \left(h(x)-h^*\right)^2 dx \\
&+\frac{qk^2}{4}\xi^2 + \frac{q}{6}w^2
\end{aligned}
$$

Inequality (63) implies the right inequality (33) with $G_1\left(s\right) := \frac{12H_{\max}}{\min\left(3H_{\max}G^{-1}\left(-cs\right), 2\mu^2, 3H_{\max}qk^2, 2H_{\max}q, 12H_{\max}g\right)}$ for $s \in [0,R)$. The proof is complete. $\triangleleft$

We next give the proof of Theorem 1.

**Proof of Theorem 1:** Let $r \in [0,R)$ be given (arbitrary) and let constants $\sigma, q, k > 0$ for which (27) holds (but otherwise arbitrary) be given.

Consider a classical solution of the PDE-ODE system (7)-(10), (12) with (28) that satisfies $V(\xi(0), w(0), h[0], v[0]) \le r$. Then we obtain from Lemma 2 and (7), (17), (28) for all $t > 0$:

$$\frac{d}{dt}V\left(\xi(t),w(t),h[t],v[t]\right)=-\mu g\int_0^L h_x^2(t,x)dx-\mu\int_0^L h(t,x)v_x^2(t,x)dx-qk^3\xi^2(t)$$
$$+qk\left(w(t)+k\xi(t)\right)^2+f(t)\left(\int_0^L\left(2h(t,x)v(t,x)+\mu h_x(t,x)\right)dx-q\left(w(t)+k\xi(t)\right)\right) \tag{64}$$

Using (28) and (64), we obtain for all $t>0$:

$$\frac{d}{dt}V\left(\xi(t),w(t),h[t],v[t]\right)=-\mu g\int_0^L h_x^2(t,x)dx$$
$$-\mu\int_0^L h(t,x)v_x^2(t,x)dx-qk^3\xi^2(t)+qk\left(w(t)+k\xi(t)\right)^2 \tag{65}$$
$$-\sigma\left(\int_0^L\left(2h(t,x)v(t,x)+\mu h_x(t,x)\right)dx-q\left(w(t)+k\xi(t)\right)\right)^2$$

Using the inequality

$$2\left(w(t)+k\xi(t)\right)\left(\int_0^L\left(2h(t,x)v(t,x)+\mu h_x(t,x)\right)dx\right)$$
$$\leq\varepsilon\left(w(t)+k\xi(t)\right)^2+\varepsilon^{-1}\left(\int_0^L\left(2h(t,x)v(t,x)+\mu h_x(t,x)\right)dx\right)^2$$

that holds for every $\varepsilon>0$, we obtain from (65) the following estimate for all $t>0$, $\varepsilon\in(0,q)$:

$$\frac{d}{dt}V\left(\xi(t),w(t),h[t],v[t]\right)\leq-\mu g\int_0^L h_x^2(t,x)dx-\mu\int_0^L h(t,x)v_x^2(t,x)dx$$
$$-qk^3\xi^2(t)-q\left(\sigma(R-\varepsilon)-k\right)\left(w(t)+k\xi(t)\right)^2 \tag{66}$$
$$+\sigma\left(\varepsilon^{-1}q-1\right)\left(\int_0^L\left(2h(t,x)v(t,x)+\mu h_x(t,x)\right)dx\right)^2$$

Using the inequality

$$2\left(2\int_0^L h(t,x)v(t,x)dx\right)\left(\mu\int_0^L h_x(t,x)dx\right)$$
$$\leq 4\zeta\left(\int_0^L h(t,x)v(t,x)dx\right)^2+\zeta^{-1}\mu^2\left(\int_0^L h_x(t,x)dx\right)^2$$

that holds for every $\zeta>0$, we obtain from (66) the following estimate for all $t>0$, $\varepsilon\in(0,q)$, $\zeta>0$:

$$
\begin{aligned}
&\frac{d}{dt}V\left(\xi(t),w(t),h[t],v[t]\right)\leq -\mu g\int_0^L h_x^2(t,x)dx-\mu\int_0^L h(t,x)v_x^2(t,x)dx \\
&-qk^3\xi^2(t)-q\left(\sigma(q-\varepsilon)-k\right)\left(w(t)+k\xi(t)\right)^2 \\
&+4\sigma\left(\varepsilon^{-1}q-1\right)\left(1+\zeta\right)\left(\int_0^L h(t,x)v(t,x)dx\right)^2 \\
&+\sigma\mu^2\left(\varepsilon^{-1}q-1\right)\left(1+\zeta^{-1}\right)\left(\int_0^L h_x(t,x)dx\right)^2
\end{aligned}
\tag{67}
$$

The Cauchy-Schwarz inequality in conjunction with (12) implies the following inequalities:

$$
\left(\int_0^L h_x(t,x)dx\right)^2\leq L\int_0^L h_x^2(t,x)dx
$$

$$
\left(\int_0^L h(t,x)v(t,x)dx\right)^2\leq m\int_0^L h(t,x)v^2(t,x)dx
$$

Thus, we obtain from (67) the following estimate for all $t>0$, $\varepsilon\in(0,q)$, $\zeta>0$:

$$
\begin{aligned}
&\frac{d}{dt}V\left(\xi(t),w(t),h[t],v[t]\right)\leq -\mu\left(g-\sigma\mu L\left(\varepsilon^{-1}q-1\right)\left(1+\zeta^{-1}\right)\right)\int_0^L h_x^2(t,x)dx \\
&-\mu\int_0^L h(t,x)v_x^2(t,x)dx-qk^3\xi^2(t)-q\left(\sigma(q-\varepsilon)-k\right)\left(w(t)+k\xi(t)\right)^2 \\
&+4m\sigma\left(\varepsilon^{-1}q-1\right)\left(1+\zeta\right)\int_0^L h(t,x)v^2(t,x)dx
\end{aligned}
\tag{68}
$$

Define for each $t\geq 0$:

$$
h_{\max}(t):=\max_{0\leq x\leq L}\left(h(t,x)\right),\ h_{\min}(t):=\min_{0\leq x\leq L}\left(h(t,x)\right)>0 \tag{69}
$$

The fact that $h_{\min}(t)>0$ follows from the fact that $h\in C^1\left([0,+\infty)\times[0,L];(0,+\infty)\right)$ (which implies that for every $t\geq 0$, $h[t]$ is a continuous positive function on $[0,L]$). Definition (69) in conjunction with (68) implies the following estimate for all $t>0$, $\varepsilon\in(0,q)$, $\zeta>0$:

$$
\begin{aligned}
&\frac{d}{dt}V\left(\xi(t),w(t),h[t],v[t]\right)\leq -\mu\left(g-\sigma\mu L\left(\varepsilon^{-1}q-1\right)\left(1+\zeta^{-1}\right)\right)\int_0^L h_x^2(t,x)dx \\
&-\mu\int_0^L h(t,x)v_x^2(t,x)dx-qk^3\xi^2(t)-q\left(\sigma(q-\varepsilon)-k\right)\left(w(t)+k\xi(t)\right)^2 \\
&+4m\sigma\left(\varepsilon^{-1}q-1\right)\left(1+\zeta\right)h_{\max}(t)\int_0^L v^2(t,x)dx
\end{aligned}
\tag{70}
$$

Since $v(t,0)=v(t,L)=0$ (recall (10)), by virtue of Wirtinger's inequality and definition (69), we have for $t\geq 0$:

$$\int_0^L v^2(t,x)dx \le \frac{L^2}{\pi^2}\int_0^L v_x^2(t,x)dx \le \frac{L^2}{\pi^2 h_{\min}(t)}\int_0^L h(t,x)v_x^2(t,x)dx \tag{71}$$

Combining (70) and (71), we obtain the following estimate for all $t>0$, $\varepsilon\in(0,q)$, $\zeta>0$:

$$\begin{aligned}\frac{d}{dt}V\left(\xi(t),w(t),h[t],v[t]\right) \le & -\mu\left(g-\sigma\mu L\left(\varepsilon^{-1}q-1\right)\left(1+\zeta^{-1}\right)\right)\int_0^L h_x^2(t,x)dx \\ & -\left(\mu-4m\sigma\left(\varepsilon^{-1}q-1\right)\left(1+\zeta\right)\frac{L^2 h_{\max}(t)}{\pi^2 h_{\min}(t)}\right)\int_0^L h(t,x)v_x^2(t,x)dx \\ & -qk^3\xi^2(t)-q\left(\sigma(q-\varepsilon)-k\right)\left(w(t)+k\xi(t)\right)^2\end{aligned} \tag{72}$$

Setting $\varepsilon = q - q\sigma^{-1}\theta\dfrac{G^{-1}(-cr)}{b+G^{-1}(-cr)}$ (notice that since $\theta=\dfrac{\sigma g}{g+\mu\sigma L}<\sigma$ it is automatically guaranteed that $\varepsilon\in(0,q)$) and $\zeta=\dfrac{\frac{\sigma\mu L}{g+\mu\sigma L}G^{-1}(-cr)}{b\varphi+(\varphi-1)\frac{\mu\sigma L}{g+\mu\sigma L}G^{-1}(-cr)}$ where $\varphi:=\dfrac{\frac{\sigma\mu L}{g+\mu\sigma L}G^{-1}(-cr)}{2b+\frac{2\mu\sigma L}{g+\mu\sigma L}G^{-1}(-cr)}+\dfrac{1}{2}<1$ (notice that $\zeta>0$), it follows from (72) and the facts that $\theta=\dfrac{\sigma g}{g+\mu\sigma L}$, $b=\dfrac{4mL^2H_{\max}}{\mu\pi^2}\theta$ that the following estimate holds for all $t>0$:

$$\begin{aligned}\frac{d}{dt}V\left(\xi(t),w(t),h[t],v[t]\right) \le & -\mu g\left(1-\varphi\right)\int_0^L h_x^2(t,x)dx \\ & -\beta\frac{\mu G^{-1}(-cr)}{H_{\max}}\left(\frac{H_{\max}}{\beta G^{-1}(-cr)}-\frac{h_{\max}(t)}{h_{\min}(t)}\right)\int_0^L h(t,x)v_x^2(t,x)dx \\ & -qk^3\xi^2(t)-q^2\left(1-\gamma\right)\theta\frac{G^{-1}(-cr)}{b+G^{-1}(-cr)}\left(w(t)+k\xi(t)\right)^2\end{aligned} \tag{73}$$

where

$$\gamma:=\frac{k\left(b+G^{-1}(-cr)\right)}{q\theta G^{-1}(-cr)} \text{ and } \beta:=\frac{b\varphi}{b\varphi+(\varphi-1)\frac{\mu\sigma L}{g+\mu\sigma L}G^{-1}(-cr)}<1 \tag{74}$$

Notice that (27) implies that $\gamma\in(0,1)$. Inequality (73) in conjunction with the facts that $\gamma\in(0,1)$, $\varphi\in(0,1)$ shows that the following implication holds:

$$\text{“If } t>0 \text{ and } \frac{h_{\max}(t)}{h_{\min}(t)}<\frac{H_{\max}}{\beta G^{-1}(-cr)} \text{ then } \frac{d}{dt}V\left(\xi(t),w(t),h[t],v[t]\right)\le 0\text{”} \tag{75}$$

The fact that $h\in C^1\left([0,+\infty)\times[0,L];(0,+\infty)\right)$ implies that both mappings $t\to h_{\max}(t)$, $t\to h_{\min}(t)>0$, defined by (69) are continuous (see Proposition 2.9 on page 21 in [14]). Consequently, the mapping $t\to\dfrac{h_{\max}(t)}{h_{\min}(t)}\ge 1$ is continuous with $\dfrac{h_{\max}(0)}{h_{\min}(0)}\le\dfrac{G^{-1}(cr)}{G^{-1}(-cr)}$. The latter inequality is a consequence of the fact that $V(\xi(0),w(0),h[0],v[0])\le r$ and (25). Since $G^{-1}(cr)<H_{\max}$ and $\beta<1$ it follows that $\frac{h_{\max}(0)}{h_{\min}(0)}<\frac{H_{\max}}{\beta G^{-1}(-cr)}$. Therefore, by continuity there exists $T>0$ such that $\frac{h_{\max}(t)}{h_{\min}(t)}<\frac{H_{\max}}{\beta G^{-1}\left(-cr\right)}$ for all $t\in[0,T)$.

We next prove by contradiction that $\frac{h_{\max}(t)}{h_{\min}(t)}<\frac{H_{\max}}{\beta G^{-1}(-cr)}$ for all $t\ge 0$. Assume the contrary, i.e. that there exists $t\ge 0$ such that $\frac{h_{\max}(t)}{h_{\min}(t)}\ge\frac{H_{\max}}{\beta G^{-1}(-cr)}$. Therefore, the set $A:=\left\{t\ge 0:\frac{h_{\max}(t)}{h_{\min}(t)}\ge\frac{H_{\max}}{\beta G^{-1}(-cr)}\right\}$ is non-empty and bounded from below. Thus we can define $t^*=\inf(A)$. By definition, it holds that $t^*\ge T>0$ and $\frac{h_{\max}(t)}{h_{\min}(t)}<\frac{H_{\max}}{\beta G^{-1}(-cr)}$ for all $t\in[0,t^*)$. By continuity of the mapping $t\to\frac{h_{\max}(t)}{h_{\min}(t)}$ we obtain that $\frac{h_{\max}(t^*)}{h_{\min}(t^*)}=\frac{H_{\max}}{\beta G^{-1}(-cr)}$. Moreover, since $\frac{h_{\max}(t)}{h_{\min}(t)}<\frac{H_{\max}}{\beta G^{-1}(-cr)}$ for all $t\in[0,t^*)$, it follows from implication (75) that $\frac{d}{dt}V\left(\xi(t),w(t),h[t],v[t]\right)\le 0$ for all $t\in(0,t^*)$. By continuity of the mapping $t\to V\left(\xi(t),w(t),h[t],v[t]\right)$, we obtain that $V\left(\xi(t^*),w(t^*),h[t^*],v[t^*]\right)\le V\left(\xi(0),w(0),h[0],v[0]\right)\le r$. On the other hand, the previous inequality and (25) imply that $\frac{h_{\max}(t^*)}{h_{\min}(t^*)}\le\frac{G^{-1}(cr)}{G^{-1}(-cr)}<\frac{H_{\max}}{\beta G^{-1}(-cr)}$ which contradicts the equation $\frac{h_{\max}(t^*)}{h_{\min}(t^*)}=\frac{H_{\max}}{\beta G^{-1}(-cr)}$.

Since $\frac{h_{\max}(t)}{h_{\min}(t)}<\frac{H_{\max}}{\beta G^{-1}(-cr)}$ for all $t\ge 0$, we conclude from implication (75) that $\frac{d}{dt}V\left(\xi(t),w(t),h[t],v[t]\right)\le 0$ for all $t>0$. By continuity of the mapping $t\to V\left(\xi(t),w(t),h[t],v[t]\right)$, we obtain that

$$V\left(\xi(t),w(t),h[t],v[t]\right)\le V\left(\xi(0),w(0),h[0],v[0]\right)\le r<R \text{ for all } t\ge 0$$

Hence, $(\xi(t),w(t),h[t],v[t])\in X$ for all $t\ge 0$ (recall definitions (16), (23)). Moreover, (25) implies that $\frac{h_{\max}(t)}{h_{\min}(t)}\le\frac{G^{-1}(cr)}{G^{-1}(-cr)}$ for all $t\ge 0$. The previous inequality in conjunction with (73) gives for all $t>0$:

$$\frac{d}{dt}V\left(\xi(t),w(t),h[t],v[t]\right)\le-\omega\left(\int_0^L h_x^2(t,x)dx+\int_0^L h(t,x)v_x^2(t,x)dx+\xi^2(t)+\left(w(t)+k\xi(t)\right)^2\right) \quad (76)$$

where

$$\omega := \min\left( \mu g(1-\varphi), \mu\left(1-\frac{\beta}{H_{\max}} G^{-1}(cr)\right), qk^3, q^2(1-\gamma)\theta \frac{G^{-1}(-cr)}{b+G^{-1}(-cr)} \right) \tag{77}$$

It follows from Lemma 3 and (76) that the following estimate holds for all $t > 0$:

$$\frac{d}{dt} V\left(\xi(t), w(t), h[t], v[t]\right) \le -\frac{\omega}{\Gamma\left(V\left(\xi(t), w(t), h[t], v[t]\right)\right)} V\left(\xi(t), w(t), h[t], v[t]\right) \tag{78}$$

where $\Gamma : [0, R) \to (0, +\infty)$ is the non-decreasing function involved in (32). Since $\Gamma : [0, R) \to (0, +\infty)$ is non-decreasing and since $V\left(\xi(t), w(t), h[t], v[t]\right) \le r$ for all $t \ge 0$, we obtain from (78) the following estimate for all $t > 0$:

$$\frac{d}{dt} V\left(\xi(t), w(t), h[t], v[t]\right) \le -\frac{\omega}{\Gamma(r)} V\left(\xi(t), w(t), h[t], v[t]\right) \tag{79}$$

By continuity of the mapping $t \to V\left(\xi(t), w(t), h[t], v[t]\right)$, the differential inequality (79) implies the following estimate for all $t \ge 0$:

$$V\left(\xi(t), w(t), h[t], v[t]\right) \le \exp\left(-\frac{\omega t}{\Gamma(r)}\right) V\left(\xi(0), w(0), h[0], v[0]\right) \tag{80}$$

Estimate (80) in conjunction with Lemma 4 implies the following estimate for all $t \ge 0$:

$$\begin{aligned} &\left\| (\xi(t), w(t), h[t] - h^* \chi_{[0,L]}, v[t]) \right\|_X^2 \exp\left(\frac{\omega t}{\Gamma(r)}\right) \\ &\le G_1\left(V\left(\xi(t), w(t), h[t], v[t]\right)\right) G_2\left(V\left(\xi(0), w(0), h[0], v[0]\right)\right) \left\| (\xi(0), w(0), h[0] - h^* \chi_{[0,L]}, v[0]) \right\|_X^2 \end{aligned} \tag{81}$$

where $G_i : [0, R) \to (0, +\infty)$, $i = 1, 2$, are the non-decreasing functions involved in (33). Since $G_i : [0, R) \to (0, +\infty)$, $i = 1, 2$, are non-decreasing functions and since $V\left(\xi(t), w(t), h[t], v[t]\right) \le r$ for all $t \ge 0$, we obtain from (81) the following estimate for all $t \ge 0$:

$$\left\| (\xi(t), w(t), h[t] - h^* \chi_{[0,L]}, v[t]) \right\|_X^2 \le G_1(r) G_2(r) \exp\left(-\frac{\omega t}{\Gamma(r)}\right) \left\| (\xi(0), w(0), h[0] - h^* \chi_{[0,L]}, v[0]) \right\|_X^2 \tag{82}$$

Estimate (29) is a direct consequence of estimate (82). The proof is complete. ⊲

We end this section by providing the proof of Proposition 1.

**Proof of Proposition 1:** Let arbitrary $(\xi, w, h, v) \in S$ that satisfies (34) with $0 < \varepsilon < \frac{\min\left(h^*, H_{\max} - h^*\right)}{\sqrt{L}}$ be given. Using definitions (17), (18), (19), the inequalities $(h(x)v(x) + \mu h_x(x))^2 \le 2h^2(x)v^2(x) + 2\mu^2 h_x^2(x)$, $\left(w + k\xi\right)^2 \le 2w^2 + 2k^2\xi^2$ and (25), we obtain:

$$V(\xi,w,h,v)=\frac{1}{2}\int_0^L h^{-1}(x)(h(x)v(x)+\mu h_x(x))^2dx+g\int_0^L\left(h(x)-h^*\right)^2dx$$
$$+\frac{1}{2}\int_0^L h(x)v^2(x)dx+\frac{qk^2}{2}\xi^2+\frac{q}{2}\left(w+k\xi\right)^2 \tag{83}$$
$$\le\frac{3}{2}\int_0^L h(x)v^2(x)dx+\mu^2\int_0^L h^{-1}(x)h_x^2(x)dx+g\int_0^L\left(h(x)-h^*\right)^2dx+\frac{3qk^2}{2}\xi^2+qw^2$$

Since $\int_0^L h(x)dx=m$ (recall definition (16)) and since $h^*=m/L$, it follows from continuity of $h$ and the mean value theorem that there exists $x^*\in[0,L]$ such that $h(x^*)=h^*$. Using the Cauchy-Schwarz inequality, we get (55) for all $x\in[0,L]$. Consequently, we obtain for all $x\in[0,L]$ from (55), (26), (34) and the fact that $0<\varepsilon<\frac{\min\left(h^*,H_{\max}-h^*\right)}{\sqrt{L}}$:

$$0<h^*-\sqrt{L}\left\|h'\right\|_2\le h(x)\le h^*+\sqrt{L}\left\|h'\right\|_2\le h^*+\varepsilon\sqrt{L}<H_{\max} \tag{84}$$

Therefore, we obtain from (83) and (84):

$$V(\xi,w,h,v)\le\frac{3H_{\max}}{2}\int_0^L v^2(x)dx+\mu^2\left(h^*-\varepsilon\sqrt{L}\right)^{-1}\int_0^L h_x^2(x)dx+g\int_0^L\left(h(x)-h^*\right)^2dx+qw^2+\frac{3qk^2}{2}\xi^2 \tag{85}$$

Inequality (38) is a direct consequence of (85) and definition (26). The proof is complete. ◁

## 6. Concluding Remarks

By applying the CLF methodology we have managed to achieve semi-global stabilization results for the viscous Saint-Venant liquid-tank system. As mentioned in the Introduction, to our knowledge, this is the first paper in the literature that achieves stabilization of the nonlinear viscous Saint-Venant system. Moreover, this is the first paper that guarantees a spill-free condition for the movement of the tank.

The obtained results leave some open problems which will be the topic of future research:
1) The results were applied to classical solutions. It is of interest to relax this to weak solutions. It is an open problem to show that the proposed feedback laws preserve their strong stabilizing role in the case of weak solutions. Moreover, it is an open problem to show existence/uniqueness of (weak) solutions for the closed-loop system. To this purpose, ideas utilized in [30] will be employed.
2) The construction of CLFs which can allow the derivation of stability estimates in stronger spatial norms for the liquid level and velocity profiles.

Another more demanding problem that can be studied in the future is the spill-free, slosh-free and smash-free movement of a glass of water. This problem arises when we want to move the glass of water to a position which is close to a wall. In this case, we need to have control on the overshoot of the glass position error in order to avoid smashing the glass on the wall.